\documentclass[smallextended]{svjour3}                        
\usepackage{graphicx}
\usepackage{amsmath,amssymb}
\usepackage{graphicx,float}
\usepackage{authblk}
\usepackage[numbers,sort&compress]{natbib}
\smartqed

\journalname{Journal of Mathematical Chemistry}

\begin{document}

\title{Analytic results on the polymerisation random graph model}

\author{Ivan Kryven}

\institute{I. Kryven \at
              University of Amsterdam, Science Park 904,  Amsterdam, The Netherlands\\
              Tel.:~+31 20 525 6423\\
              Fax.:+31 20 525 5604\\
              \email{i.kryven@uva.nl}          
}


\maketitle

\begin{abstract}
The step-growth polymerisation of a mixture of arbitrary-functional monomers is viewed as a time-continuos random graph process with degree bounds that are not necessarily the same for different vertices. The sequence of degree bounds acts as the only input parameter of the model. This parameter entirely defines the timing of the phase transition.
Moreover, the size distribution of connected components features a rich temporal dynamics that includes: switching between exponential and algebraic asymptotes and acquiring oscillations.
The results regarding the phase transition and the expected size of a connected component are obtained in a closed form. An exact expression for the size distribution is resolved up to the convolution power and is computable in subquadratic time.
The theoretical results are illustrated on a few special cases, including a comparison with Monte Carlo simulations.\\

\end{abstract}

\section{Introduction}
The chemical graph theory is the branch of mathematical chemistry that applies graph theory to mathematical modelling of chemical processes.  This theory centres its attention on the concept of a molecular graph, which identifies atoms (or monomers) as vertices and chemical bonds as edges. This structure, finite or infinite, is usually defined a priori, \emph{e.g.} molecular graphs  describing structural isomers or Euclidian graphs describing crystal nets \cite{bonchev1991,rouvray2002}. The graph-theoretical invariants of such chemical objects are known to be strongly correlated with physical properties of the resulting materials. These invariants include but are not restricted to: Wiener index, average shortest path, shape index, centric index, and connectivity index \cite{mihalic1992,mohar1988,randic1991,estrada2001}.
Not all molecular topologies can be described by a single graph, but rather by a probability measure over  graphs \cite{kryven2015c,kryven2013c}. This scope covers (hyper-)branched polymers, cross-linked polymers, molecular networks, and gels to name a few. A branch of graph theory that operates with probability distributions over graphs -- random graph theory -- has little documented applications to chemistry at present. 

Consider a chemical system where each monomer has a predefined functionality, that is the maximum number of neighbours in the network. If the spatial positioning of the monomers  is disregarded, the monomers can be represented as vertices in a graph model. From this perspective, the polymerisation process is a random graph process that respects the limitations induced by the chemistry, for instance, the bound on the vertex degree. The fact that this chemical system can be well described by graph theory is already hinted by a broad range of  analogues to graph-theoretical terminology that exists in polymer chemistry: vertex (\emph{monomer}), degree bound (\emph{functionality}), graph (\emph{polymer network}), tree (\emph{branched polymer}), connected component (\emph{polymer molecule}), giant component (\emph{gel}), density (\emph{conversion}), etc.

In this paper a random graph process is introduced to model an evolving molecular network. The degree distribution of this random graph is defined by a time-continuous evolution equation that mimics the chemistry of the step-growth polymerisation process. This process starts with disconnected vertices and progresses up to the point where no new edge can be placed.
The degree of each vertex is bounded, but different bounds may be defined for distinct vertices.
Therefore, we distinguish between the degree -- actual number of incident edges, and the functionality -- pre-imposed bound on the number of incident edges.
At each time step, the probability that a vertex receives an edge is proportional to the difference between  the vertex's functionality and degree.
The share of vertices in each functionality class is pre-defined, and constitutes the only input parameter for the random graph model. 

Most of the available studies target narrow special cases of this system and pursue results with a distinct reasoning from the graph-theoretical one. Important contributions include:
Hamilton-Jacobi formalism as applied to dynamic graphs with globally bounded degrees \cite{ben2011}, results on the grabbing-particle system \cite{bertoin2010}, open-form analytical results for non-phase-transiting systems \cite{hillegers2015}, combinatorial analysis for monomers bearing identical groups \cite{durand1982}, closed-form analytical \cite{zhou2006} and numerical \cite{kryven2013c} results for trifunctional vertices in a directed topology,  analytical results for mixture of bi- and trifunctional vertices \cite{iedema2012a}, analytical results on phase transition in evolving directed graphs \cite{kryven2016}, and stochastic simulations on molecular networks \cite{kryven2016b}. The random graph model is also related to many processes outside polymer chemistry. For instance, Smoluchowski coagulation equation with a multiplicative kernel governs the dynamics of component-size distribution of the polymerisation random graph with trifunctional vertices. Only in this special case, the analytical expression for component sizes is available also after the phase transition, for a review on  Smoluchowski coagulation see Refs. \cite{wattis2006,bertoin2009}. In probability theory, the  gambler's ruin problem for infinite number of games is equivalent to finding criteria for the phase transition in the polymerisation random graph with vertices not exceeding degree three \cite{harik1999}.

The rest of the  paper is organised as follows. First, a differential-difference equation describing evolution of the degree distribution due to the step-growth polymerisation process is formulated and solved in time. Then, given the time dependent degree distribution, the emergence of the giant component is analysed. This includes results on the edge density at which the giant component appears and the criterion on the functionality distribution that admit emergence of the giant component at finite time. Furthermore, the size distribution of connected components is resolved  and expressions for the expected component size are given. Finally, the theoretical results are discussed for a few special cases. The theory is also compared against the size distributions that were generated by a Monte Carlo simulation.

\section{Evolution process for the degree distribution}
This paper studies infinite graphs as a model for a polymer network: a chemical system composed of randomly interconnected identical units. 
In the infinite graph, degree distribution $\text{u}(n), n=0,1,2,\dots$ is the probability that a randomly sampled vertex has $n$ adjacent edges~\cite{newman2010}. Since a degree of a vertex cannot be arbitrary large in a chemical system, each vertex is assigned a bound on its degree, $m=0,1,2,\dots$. To copy the chemical terminology, we refer to this bound as the functionality~\cite{Stockmayer1943}.  So that one may speak of a two-variate distribution $\text{u}(n,m), \;n,m=0,1,2,\dots$ as the probability to sample a vertex with degree $n$ and functionality $m$, such that $\text{u}(n,m)=0$ for $n>m$. 
We will now construct an evolutionary process for $\text{u}(k)$ that mimics the step-growth polymerisation of multifunctional monomers. 
This linking process starts with disconnected vertexes, that is the probability to sample a vertex  of degree zero is $d(0,k)=1,$ and the process ends when one samples a vertex with $n=m$ with probability one. The precise rule of assigning a new edge is the following conceptualisation of the step-growth polymerisation process: on each time step, one samples two candidate vertices with probability proportional to $(m-n)\text{u}(n,m)$ and connects them with an edge. So that
\begin{equation}\label{eq:process1}
\{(n_1,m_1),(n_2,m_2)\}\rightarrow\{(n_1+1,m_1),(n_2+1,m_2)\},\; n_1\leq m_1,\; n_2\leq m_2,
\end{equation}
where  $(n_1,m_1)$ and  $(n_2,m_2)$ are the configurations of the candidate vertices.
This linking process may be viewed as a generalisation of the linking process with constant degree bounds (all vertices have the same functionality $m$) as introduced in Ref. \cite{ben2011}, Eq. (3). An alternative way of introducing \eqref{eq:process1} is by writing the corresponding reaction mechanism  for monomer species $M_{n,m}$:
\begin{equation}\label{eq:mechanism}
 M_{n_1,m_1} + M_{n_2,m_2}  \xrightarrow{ (m_1 - n_1 ) ( m_2 - n_2) } M_{ n_1 + 1, m_1} + M_{n_2+1,m_2}.
\end{equation}
Both notations \eqref{eq:process1} and \eqref{eq:mechanism} are equivalent and correspond to the following Kolmogorov forward equation governing the  evolution of $\text{u}(n,m)$,
\begin{equation}\label{eq:population.balance}
\begin{aligned}
\frac{\partial}{\partial t} \text{u}(n,m,t) = &  \Big( (m - n + 1) \text{u}(n - 1,m,t) - ( m - n) \text{u}(n,m,t) \Big)\times\\
  &  \sum\limits_{m=0}^\infty\sum\limits_{n=0}^m \Big(m\, \text{u}(n,m,t) -  n\, \text{u}(n,m,t)\Big);\\
\end{aligned}
\end{equation}
where at $t=0$, $\text{u}(n,m,t)$ satisfies the following initial conditions,
\begin{equation}\label{eq:init}
\begin{aligned}
 \text{u}(0,m,0)&= f_m,\\
 \text{u}(n,m,0)&= 0,\;n>0.
 \end{aligned}
\end{equation}
 In this equation, the probability to sample a vertex of functionality $m$ is constant over time, $\sum\limits_{n=0}^{\infty} u(n,m,t) = f_m,\; \sum\limits_{m=1}^\infty f_m=1$,
 and $f_m$ is treated as the only parameter of the model.
 The sum written in the second line of Eq. \eqref{eq:population.balance} represents the expected number of unused but potentially available edges and can be viewed as a difference of two partial moments, $\mu_{01}(t)-\mu_{10}(t),$ where   
\begin{equation}\label{eq:mom}
\begin{aligned}
\mu_{ij}(t) =&\sum\limits_{m=0}^\infty\sum\limits_{n=0}^m n^i \,m^j\, \text{u}(n,m,t). \\
\end{aligned}
\end{equation}
The edge density, $c(t)\in [0, 1],$ is a ratio of expected number of edges at time $t$ to the expected number of edges at the end of the process:
$$
c(t)=\frac{\mu_{10}(t) }{\mu_{01} }.
$$ 
It is convenient to use $c(t)$ as an alternative measure of the progress. 
The differential equation \eqref{eq:population.balance} falls into the class of linear population balance equations.  This class of equations  frequently appears as a model for many chemical and biological problems where it is usually approached numerically \cite{kryven2015a,kryven2014d}.
In the current case, it is possible to find an analytical solution of \eqref{eq:population.balance} by transforming the equation to the domain of generating functions, solving the corresponding partial differential equation, and applying the inverse transform. 

Let us rewrite \eqref{eq:population.balance} in terms of a univariate generating function, $$\text{U}(x,m,t)=\sum\limits_{n=0}^m x^n \text{u}(n,m,t), \; |x|<1, \; x\in\mathbb{C}.$$
Taking generating function transform on both sides of Eq. \eqref{eq:population.balance} leads to a partial differential equation (PDE),
\begin{equation}\label{eq:pde}
\left\{
\begin{aligned}
\frac{\partial}{\partial t} \text{U}(x,m,t) =& \Big( (m x - m) \text{U}(x,m,t) + (x - x^2) \frac{\partial}{\partial x} \text{U}(x,m,t) \Big) (\mu_{01}(t) - \mu_{10}(t))\\
\text{U}(x,m,0) =& f_m,\; |x|<1.
\end{aligned}\right.
\end{equation}
The first partial moments appearing in \eqref{eq:pde} can be related to the generating function $\text{U}(x,m,t)$,
\begin{equation}\label{eq:genmom}
\begin{aligned}
\mu_{01}(t) = & \sum\limits_{m=0}^\infty m \text{U}(1,m,t); \\
\mu_{10}(t) = & \sum\limits_{m=0}^\infty \frac{\partial}{\partial x} \text{U}(x,m,t)|_{x=1}; \\
\end{aligned}
\end{equation}
Substituting \eqref{eq:genmom} into \eqref{eq:pde} we obtain a system of ordinary differential equations for the partial moments,
\begin{equation}\label{eq:ode}
\left\{
\begin{aligned}
\mu'_{10}(t) =& \big(\mu_{01}(t) - \mu_{10}(t)\big)^2,\\
\mu'_{01}(t) = &0, 
\end{aligned}\right.
\end{equation}
that is subject to initial conditions $\mu_{10}(0) = 0, \;\mu_{01}(0) = \mu_{01}$.
Solving \eqref{eq:ode} gives
\begin{equation}\label{eq:mom_sol}
\begin{aligned}
\mu_{10}(t) = & \frac{\mu_{01}^2 t}{1 + \mu_{01} t};\\
\mu_{01}(t) = & \mu_{01};\\
\end{aligned}
\end{equation}
Now, having  explicit expressions for $\mu_{10}(t),\mu_{01}(t)$ at hand, allows us to write the solution of PDE \eqref{eq:pde},
\begin{equation}\label{eq:unmt}
\text{U}(x,m,t)= \Big(\frac{1 + \mu_{01} t }{1 + \mu_{01} t \,x}\Big)^{-m} f_m,
\end{equation}
which, in turn, generates $\text{u}(n,m,t),$
\begin{equation}\label{eq:sol_t}
\text{u}(n,m,t) =  \binom{m}{n} (\mu_{01} t)^n (1+\mu_{01} t)^{-m}f_m,\; n\leq m.
\end{equation}
The latter expression can be reformulated in terms of edge density $c(t)$ instead of time.
To do this, it is enough to realise that $c(t)=\frac{\mu_{10}(t)}{\mu_{01}} = \frac{\mu_{01} t}{1 + \mu_{01} t}$ and 
$(\mu_{01} t)^n (1+\mu_{01} t)^{-m}=(\mu_{01} t)^n(1+\mu_{01} t)^{-n}  (1+\mu_{01} t)^n (1+\mu_{01} t)^{-m}=
\left(\frac{\mu_{01} t }{1+\mu_{01} t}\right)^n  \left(\frac{1}{1+\mu_{01} t}\right)^{m-n}=
\left(\frac{\mu_{01} t }{1+\mu_{01} t}\right)^n  \left(1-\frac{\mu_{01} t}{1+\mu_{01} t}\right)^{m-n}$ so that  Eq. \eqref{eq:sol_t} transforms to 
\begin{equation}\label{eq:binom}
\text{u}(n,m,t) = \binom{m}{n} c^n(t) \big(1-c(t)\big)^{m-n}f_m,\; n\leq m.
\end{equation}
Expressions \eqref{eq:sol_t},\eqref{eq:binom} satisfy the initial conditions \eqref{eq:init}, 
whereas in the limiting case of $t \to \infty$, the degree distribution and the distribution of maximal functionalities coincide:
$$
\begin{cases}
\lim\limits_{t\rightarrow\infty} \text{u}(n,m,t)=f_m, & n=m;\\
\lim\limits_{t\rightarrow\infty} \text{u}(n,m,t)=0, & n<m.
\end{cases}
$$
The actual degree distribution $\text{u}(n)$, is found by summating $\text{u}(n,m,t)$ over   functionalities $m,$ 
\begin{equation}\label{eq:the.solution}
\text{u}(n,t)=\sum\limits_{m=1}^\infty \text{u}(n,m,t).
\end{equation}
Here, we employed the fact, that $\text{u}(n,m,t)=0,$ for $n>m.$
Degree distribution $u(n,t)$ evolves form the Kronecker's delta function, $\delta_n$ at $t=0$ to $f_m$ in the limit of $t\rightarrow \infty.$
The moments of the degree distribution, $\mu_i=\sum\limits_{n=0}^\infty n^i u(n,t) = \mu_{i0}$ can be directly found from summation of Eq. \eqref{eq:the.solution}. For instance the expressions for the first three moments read,
\begin{equation}\label{eq:moms}
\begin{aligned}
\mu_{1}(t) &= \frac{\mu_{01}^2 t}{1 + \mu_{01} t},\\
\mu_{2}(t) &= \frac{\mu_{01}^2 t(1 +  \mu_{02} t)}{(1 + \mu_{01} t)^2},\\
\mu_{3}(t) &= \frac{ \mu_{01}^2 t (1 - 3 \mu_{01} t + 4 \mu_{02}  t) }{1 + \mu_{01} t^2}.
\end{aligned}
\end{equation}

\section{Global properties of the network, the giant component}\label{sec:giant}
Up to this point we have discussed only local properties, i.e. the way the graph can be seen from a viewpoint of a single vertex. However, in a 
randomly interconnected system, local properties, as for instance the degree distribution, play a decisive role in defining the global properties of the graphs itself. An important finding  that allows us to connect the the two worlds is the result by Molloy and Reed on the existence of the giant component \cite{molloy1998}: there exists a component of the same order of magnitude as the whole graph (the giant component) iff, 
$$
\sum\limits_{n=1}^\infty n(n-2) 	\text{u}(n,t)> 0,
$$
while the equality is reached exactly at the phase transition point. This phase transition condition can be rewritten in terms of moments \eqref{eq:genmom},
\begin{equation} \label{eq:gelcondition}
	\mu_{20}(t)-2\mu_{10}(t) = 0.
\end{equation}
Substituting the analytical expression for moments \eqref{eq:mom_sol} into Eq. \eqref{eq:gelcondition} we obtain the phase transition time (or the gelation time in the chemical terminology),
\begin{equation}\label{eq:tgel}
t_{g} = \frac{1}{\mu_{02}-2\mu_{01}}.
\end{equation}
Similarly, the edge density at the phase transition (or gel conversion) is written out as
\begin{equation}\label{eq:cgel}
c_{g} =  \frac{\mu_{01} t_g}{1 + \mu_{01} t_g}= \frac{\mu_{01}}{ \mu_{02}- \mu_{01}}. 
\end{equation}
From the last relation \eqref{eq:cgel} we can see that the system features the phase transition in a finite time only when $\mu_{02}-2\mu_{01} > 0.$ If the inequality is replaced by an equality ($\mu_{02}=2\mu_{01}$), then the phase transition  will be approached asymptotically at $t \rightarrow \infty.$
This brings us to the following, especially important for its chemical context,\\
\emph{Corollary:} let $M$ monomer species of functionalities $m=1,\dots,M$  and fractions $f_1,f_2,\dots,f_M,\; \sum\limits_{m=1}^M f_m=1$ react at constant rate $k_p$, then the system features the phase transition in a finite time if and only if 
\begin{equation}\label{eq:col1}
\sum\limits_{m=1}^M m^2 f_m -2\sum\limits_{m=1}^M m f_m>0.
\end{equation}
If the phase transition occurs, then it occurs at the following time and edge density,  
$$t_{g} = \Big( k_p \sum\limits_{m=1}^M (m^2-2m) f_m \Big)^{-1},$$
and 
\begin{equation}\label{eq:col2}
c_{g} = \Big(  \sum\limits_{m=1}^M (m^2-m) f_m \Big)^{-1}\sum\limits_{m=1}^M m f_m.
\end{equation}
As special cases of this corollary, the following statements hold true.\\
1. If all monomers have the same functionality $m$, then the phase transition is reached in a finite time only if $m\geq3$ (i.e. $m$ is the smallest positive integer satisfying $m^2-2m>0$). \\
2. Adding (or removing) monomers of functionality two does not affect phase transition time $t_g$, whereas it does alter the edge density at the phase transition, $c_g$.  \\
3. Adding sufficient amount of $f_1$ to any system will prevent the phase transition.\\
4. Consider a system that consists of two species: monomers with functionality $m$ that are present at fraction $f_m$ and monomers with functionality one, that are present at fraction $f_1=1-f_m$. 
   The system does not go through the phase transition in finite time if,  
\begin{equation}\label{eq:postpone}
	f_1 > \frac{m^2-2 m}{m^2-2m+1}. 
\end{equation}
5. When all monomers have functionality $m$, the polymerisation leads to an infinite network at edge density $$c_g=\frac{m} {m^2-m}=\frac{1}{m-1}.$$\\

The latter equation was derived by Flory \cite{flory1941}. Although Flory did not consider non-constant functionality, somewhat later, he conjectured that the equation can be generalised for a mixture of arbitrary functional monomers if $m-1$ were replaced ``by the appropriate average, weighted according to the numbers of functional groups.'' (see \cite{Flory1953},~p.~353).

\section{Size distribution of connected components}
For the sake of brevity we drop time argument $t$ where it leads to no confusion, and refer to the degree distribution, as given in Eq. \eqref{eq:the.solution}, by simply $\text{u}(n)$ or by its generating function, 
\begin{equation}\label{eq:generating}
\text{U}(x)=\sum\limits_n x^n \text{u}(n), \;|x|\leq1, \; x\in\mathbb{C}.
\end{equation}
We will now apply the theory from Refs. \cite{newman2010,newman2001} to recover other non-local properties of the polymer network.

When talking about a property of a randomly sampled vertex in an infinite graph it is important to specify what is exactly the sampling rule. Up to this point, we considered the case when every vertex has equal chances to be sampled.
Consider a different strategy to choose a vertex: suppose one samples \emph{an edge} at random, so that every edge has equal probability to be sampled.  Then, one of incident to this edge vertices is chosen and the edge itself is removed. We will refer to this vertex as the biased vertex. Let $\text{u}_1(n)$ denotes the probability that a biased vertex has $n$ incident edges. Then,  
$$\text{u}_1(n) = \frac{(n+1) \text{u}(n+1)}{\sum\limits_{n=1}^{\infty} n \text{u}(n) },$$
and the corresponding generating function is 
\begin{equation}\label{eq:U1}
\text{U}_1(x)=\frac{\text{U}'(x)}{\text{U}'(x)|_{x=1}}.
\end{equation}

A connected component is a subset of vertices in a graph, such that every couple of vertices is connected with a path. Let $w(n)$ denotes the probability that a randomly sampled node belongs to a connected component of size $n.$ Similarly to definition of $\text{u}_1(n)$, let $w_1(n)$ denotes the probability that a biased vertex  belongs to a connected component of size $n$.
Newman et al. \cite{newman2010} noticed that the generating functions for $u_1(n)$ and $w_1(n)$ are related by a functional equation
\begin{equation}\label{eq:rec}
W_1(x) = x \sum\limits_{n=0}^\infty \text{u}_1(n) W_1^n(x),
\end{equation}
where $W_1(x)$ generates $w(n)$ and $U_1(x)$ generates $u_1(n)$.
This equation has a straightforward interpretation: the equation unfolds the generating function for $w_1(n)$ as a sum over all configurations of a biased vertex. Each configuration occurs with probability  $\text{u}_1(n)$ and involves $n$ biased sub-components of size $w_1(n).$ Furthermore, the sum in Eq. \eqref{eq:rec} can be in itself viewed as the definition of the generating function. So that one may write,
\begin{equation}\label{eq:W1_1}
	W_1(x) = x U_1\Big( W_1(x) \Big).
\end{equation}
Following a similar logic to derivation of \eqref{eq:W1_1}, the generating function for $w(n)$ reads
\begin{equation}\label{eq:W_1}
	W(x) = x U\Big( W_1(x) \Big).
\end{equation}
Due to Lagrange inversion principle\cite{bergeron1998}, the system of functional equations \eqref{eq:W1_1},\eqref{eq:W_1} has a unique solution. Furthermore, the formal expression for $w(n)$ can be written out in terms of convolution powers\cite{kryven2017a},
\begin{equation}
\label{eq:Lagrange1d}
w(n,t)=\begin{cases}
\frac{\mu_{01}^2 t}{(1 + \mu_{01} t)(n-1)}u_1^{*n}(n-2), & n>1, \\
u(0) & n=1.
\end{cases}
\end{equation}
Here $u_1^{*n}(n)$ denotes the convolution power,
 \begin{equation}\label{eq:conpower}
 u(k)^{*n}=u(k)^{*n-1}*u(k),
 \end{equation}
 where
 $$f(k)*g(k) =\sum\limits_{i+j=k}f(i)g(j),\;i,j,k\geq0.$$ 
On practice, the exact numerical values of \eqref{eq:Lagrange1d} can be computed by making use if the convolution theorem and evaluating \eqref{eq:conpower} with the fast Fourier transform algorithm. Such numerical routine results in $\mathcal O( n \log n)$ multiplicative operations. If all vertices have the same functionality $m$, $f_m=1$, then $w(n)$ is simply given by,
\begin{equation}
w(n,t)=\frac{\mu_{01}^2 t}{(1 + \mu_{01} t)(n-1)} \binom{ n (m-1)}{n-2}(1 + \mu_{01} t)^{- n (m  - 3)}   (\mu_{01} t)^{n-2}.
\end{equation}

The restrictions imposed by chemistry of the polymerisation system guarantee that $u(n)=0$ for some $n>n_{max}.$ 
This class of degree distributions features a defined asymptotic behaviour of $w(n)$ at large $n\gg1,$ see Ref. \cite{kryven2017a}. 
Namely,  
$$\lim\limits_{n\to \infty}\frac{w(n)}{w_\infty(n)}=1,$$
 where
\begin{equation}\label{eq:asymptote}
w_\infty(n) = C_1(t)e^{-C_2(t) n} n^{-3/2},
\end{equation}
and the coefficients are given by
$$
\begin{aligned}
C_1(t) &=\frac{\mu_{1}^2(t)}{ \sqrt{2 \pi \Big(\mu_{1}(t) \mu_{3}(t)-\mu_{2}^2(t)\Big) }}= \frac{\mu_{01}^2 \sqrt{t}}{\sqrt{2 \pi( \mu_{02}-\mu_{01} ) (2 + 3 \mu_{01} t - \mu_{02} t)}},\\
C_2(t) &=  \frac{\Big(\mu_{2}(t)-2 \mu_{1}(t)\Big)^2}{2\Big( \mu_{1}(t) \mu_{3}(t)- \mu_{2}^2(t)\Big)}=\frac{(1 -( \mu_{02} - 2 \mu_{01})t )^2}{2 t (\mu_{02} - \mu_{01})  (2 + 3 \mu_{01} t - \mu_{02} t)}.
\end{aligned}
$$
In the latter transformation we made use of the expressions of the moments \eqref{eq:moms}.

One may see that at the phase transition, when $t=t_g$, the coefficient in the exponential function in \eqref{eq:asymptote} vanishes and the asymptote switches to the power law decay.

\section{Expected size of connected components}
It is important to note, that $w(n)$ describes only finite components. Before the phase transition, a randomly sampled node belongs to a finite component with probability one, therefore $\sum\limits_{n=1}w(n)=W(1)=1$. After the phase transition, when $t>t_{\text{g}},$  the probability that a randomly sampled node belongs to a finite component is smaller than one and
$w(n)$ fails to be normalised: 
$$W(1)=\sum\limits_n w(n)=1-g_f,$$
where  $g_f$ is the the probability  that a randomly sampled vertex belongs to the giant component (or gel fraction): $g_f=0$ for $t<t_g$ and $g_f\in[0,1]$ for $t>t_g$. 
Plugging $x=1$ into \eqref{eq:W1_1} one obtains,
\begin{equation}\label{eq:gf}
g_f = 1 - W(1) = 1 - \text{U}( r_0 ),
\end{equation} 
where $r_0:=W_1(1)$ is the smallest positive fixed point of $\text{U}_1(x),$
\begin{equation}\label{eq:r0}
r_0=\text{U}_1(r_0).
\end{equation}
We will now derive the expression for the expected size of connected component, as given by 
$$M_w:=\frac{\sum\limits_{n=1}n w(n)}{\sum\limits_{n=1} w(n)}=\frac{W'(1)}{W(1)}.$$
 Let $t<t_g$, then $W(1)=W_1(1)=1$ and evaluating $W_1'(1)$ from Eq. \eqref{eq:W1_1} gives,
$$W_1'(1)= \text{U}(W_1(1)) +\text{U}'(W_1(1))W_1'(1)=1 + \text{U}_1'(1) W_1'(1) = \frac{1}{1-\text{U}_1'(1)}.$$
Similarly, evaluating $W'(1)$ from Eq. \eqref{eq:W_1} gives,
\begin{equation}\label{eq:mw1}
\begin{aligned}
	 M_w=&W'(1)= 1 + \text{U}'(1) W_1'(1)=1 + \frac{\text{U}'( 1)}{1-\text{U}_1'(1)}= 1-\frac{\mu_{1}^2(t)}{\mu_{2}(t)-2\mu_{1}(t)} \\
	=&1+\frac{ \mu_{01}^2 t }{1 + 2 \mu_{01} t - \mu_{02} t},\;t<t_g.
	 \end{aligned}
\end{equation}
The latter transformation is made realising that $\text{U}'( 1)=\mu_{1}(t),$  $\text{U}_1'( 1)=(\mu_{2}(t)-\mu_{1}(t))/\mu_{1}(t)$ and the moments of the degree distribution are as defined by Eqs. \eqref{eq:moms}.

Let $t>t_g$, then $W_1(1)=r_0\neq1$ and evaluating $W_1'(1)$ from Eq. \eqref{eq:W1_1} gives the following equality,
 \begin{equation}
 \begin{aligned}
W_1'(1)=  \text{U}_1(W_1(1)) + \text{U}_1'(W_1(1)) W_1'(1) = r_0 + \text{U}_1'(r_0) W_1'(1),
 \end{aligned}
  \end{equation}
 so that
 $$W_1'(1)=\frac{r_0}{1-\text{U}'_1(r_0)}.$$
Evaluating  $W'(1)$ from Eq. \eqref{eq:W_1}, gives
$$W'(1)= \text{U}(W_1(1)) + \text{U}'(W_1(1)) W_1'(1) = W(1)+ \text{U}'(r_0) \frac{r_0}{1-\text{U}'_1(r_0)}.$$
Now, realising that according to Eq. \eqref{eq:U1}, $ \text{U}'(r_0) = \text{U}'(1)\text{U}_1(r_0)= \mu_{1}r_0,$ one obtains:
\begin{equation}
\begin{aligned}\label{eq:mw2}
M_w=&\frac{W'(1)}{W(1)} = \frac{1-g_f +  \frac{\mu_{1}r_0^2}{1-\text{U}'_1(r_0)}}{1-g_f}=1+\frac{\mu_{1}r_0^2}{(1-g_f)(1-\text{U}'_1(r_0))}=\\
 =&1+\frac{\mu_{01}^2r_0^2t}{(1+\mu_{01}t)(1-g_f)(1-\text{U}'_1(r_0))},\;t>t_g.
\end{aligned}
\end{equation}
Together, Eqs. \eqref{eq:mw1} and \eqref{eq:mw2} define the expected component size before and after the phase transition, that is at $t\in[0,t_g)\cup(t_g,\infty)$. Precisely at the phase transition, $t=t_g,$ the expected component size diverges, as $(t-t_g)^{-1}$. So that 
$$\lim\limits_{t\to t_g} \frac{M_w(t)}{ (t-t_g)^{-1}} =\mathcal O(1).$$
This happens due to a different type of the asymptotical behaviour of the size distribution at the phase transition, see Eq. \eqref{eq:asymptote}. 

\section{Interpretation of the results \& examples}
\begin{figure}[htbp]
\begin{center}
 \includegraphics[width=0.7\textwidth]{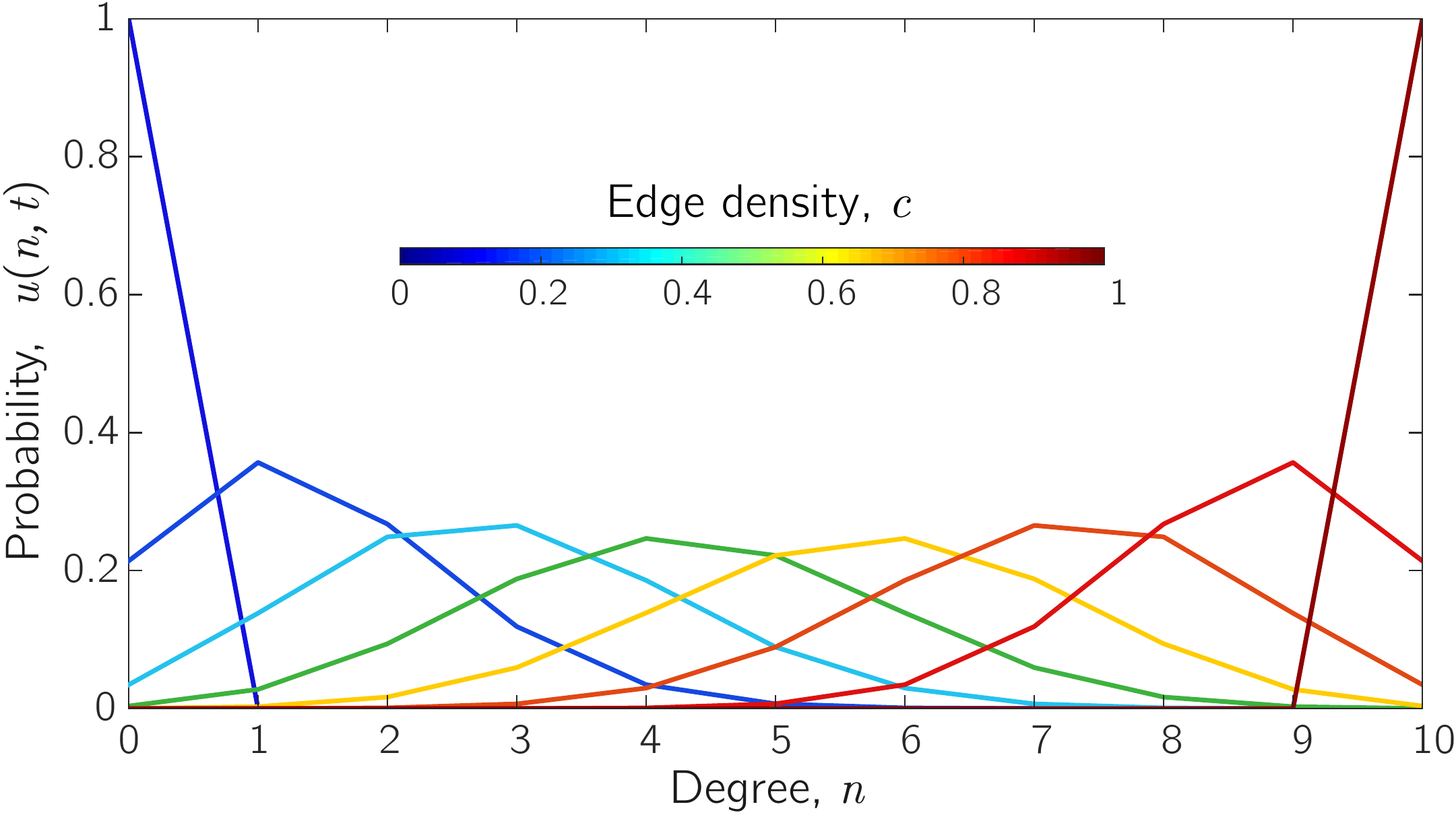} 
 \caption{Evolution of the degree distribution for $f_{10}=1$. }
\label{fig:degrees}
\end{center}
\end{figure}

The present paper introduces a model for studying polymer networks composed of multifunctional monomers that polymerise according to the step-growth mechanism \eqref{eq:mechanism}. This model associates a vertex with a monomer and an edge with a chemical bond between two such monomers in the network. A resulting topology of the polymer network is viewed as a random graph defined by its degree distribution. Initial fractions of monomers of different functionalities $f_m$ are directly related to molar concentrations of monomer species. 
The reaction kinetics is formalised by the master equation \eqref{eq:population.balance} and yields an analytical expression for the degree distribution at any point of time \eqref{eq:sol_t}. Although the master equation \eqref{eq:population.balance} has a unit rate, an arbitrary reaction rate can be modelled by simply scaling time variable $t$ in a linear fashion. An example of a degree distribution evolving in time is given in Figure~\ref{fig:degrees}. In this example, the initial condition of the kinetic model is chosen to be $f_{10}=1,$ that corresponds to pure 10-functional monomers.  In the given context, both, initial and terminal degree distributions  are Kronecker's delta functions positioned correspondingly at $m=0$ and $m=10$.

\begin{figure}[htbp]
\begin{center}
 \includegraphics[width=0.9\textwidth]{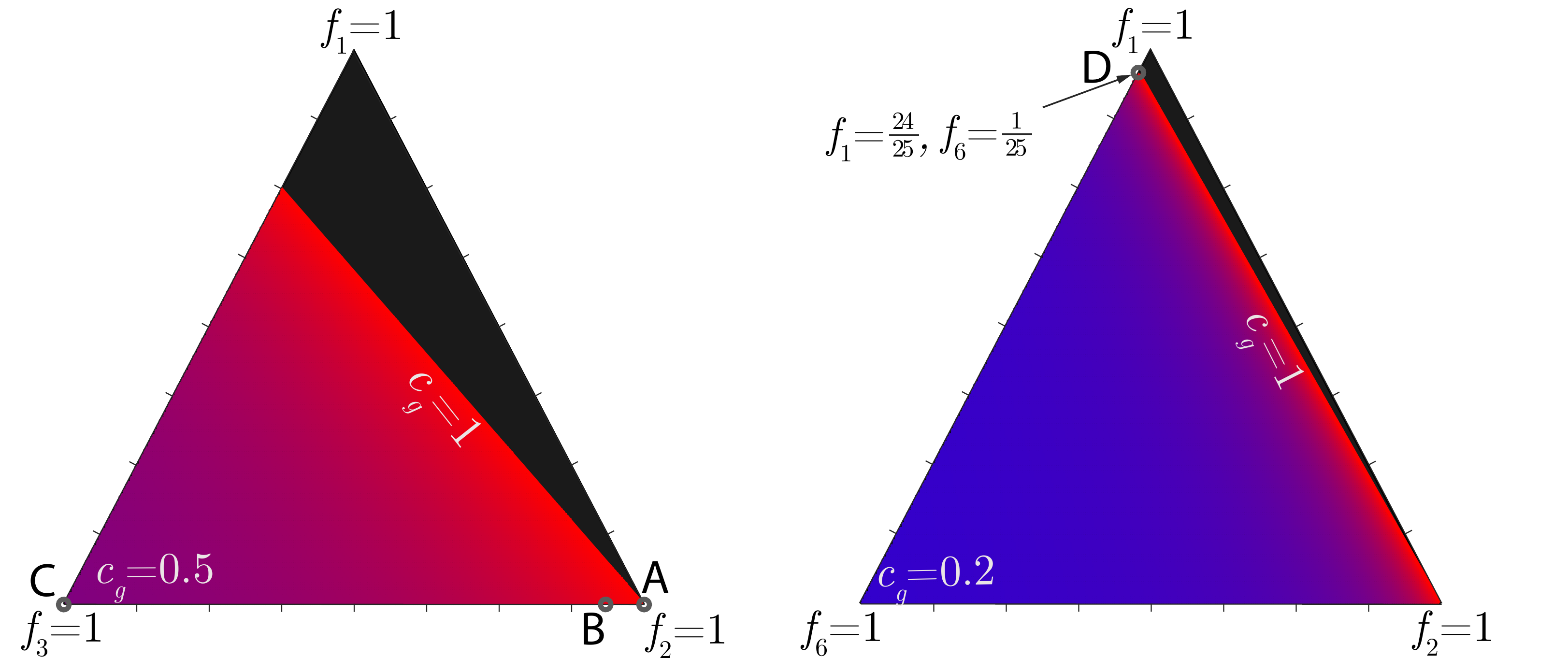} 
\caption{The edge density at phase transition, $c_g,$ is plotted as a function of concentration in barycentric coordinates for two sets of monomer functionalities: (\emph{left:}) the only non-zero concentrations are $f_1,f_2,f_3$, (\emph{right:}) the only non-zero concentrations are $f_1,f_2,f_6$. The black area corresponds to the configurations that does not feature the phase transition. The points ($A,B,C,D$) refer to special cases discussed in the paper. }
\label{fig:bar}
\end{center}
\end{figure}

A deeper analysis reveals that when initial concentrations of monomers satisfy condition \eqref{eq:col1}, the random graph develops a giant component at time $t_g$ that is given by Eq. \eqref{eq:tgel}. This event is related to the fact that the molecular network undergoes a phase transition. Such phase transition is called gelation, and is a well-documented chemical phenomenon that signifies transition from liquid-like to solid-like state in soft matter\cite{ziff1980,winter1997}. Figure~\ref{fig:bar} presents two examples showing how $t_g$ is influenced by varying $f_m.$
The figure illustrates the fact that addition of one- and two- functional vertices may be used to control the timing of the phase transition: addition of two-functional vertices postpones the emergence of the giant component in terms of $c_g,$ whereas $t_g$ remains invariant; addition of one-functional vertices may entirely prevent it.

The size distribution of connected components, as given in Eq. \ref{eq:Lagrange1d}, is interpreted as the molecular weight distribution, whereas the asymptote \eqref{eq:asymptote} might serve as a good way to approximate the latter if rapid computations are required. Evolution of the expected number of this distribution, also known in the chemical literature as number-average molecular weight, is given by Eqs. \eqref{eq:mw1},\eqref{eq:mw2}.
 
More examples of phase transitioning systems, as obtained for a few instances of functionality distribution $f_m,$ follow below.  
These examples are supplemented with a MATALB code that reproduces the size distribution and the corresponding expected value for an arbitrary functionality distribution and the process time \cite{kryven.git.polyrandgrpah}. 
 
\begin{figure}[htbp]
\begin{center}
 \includegraphics[width=0.65\textwidth]{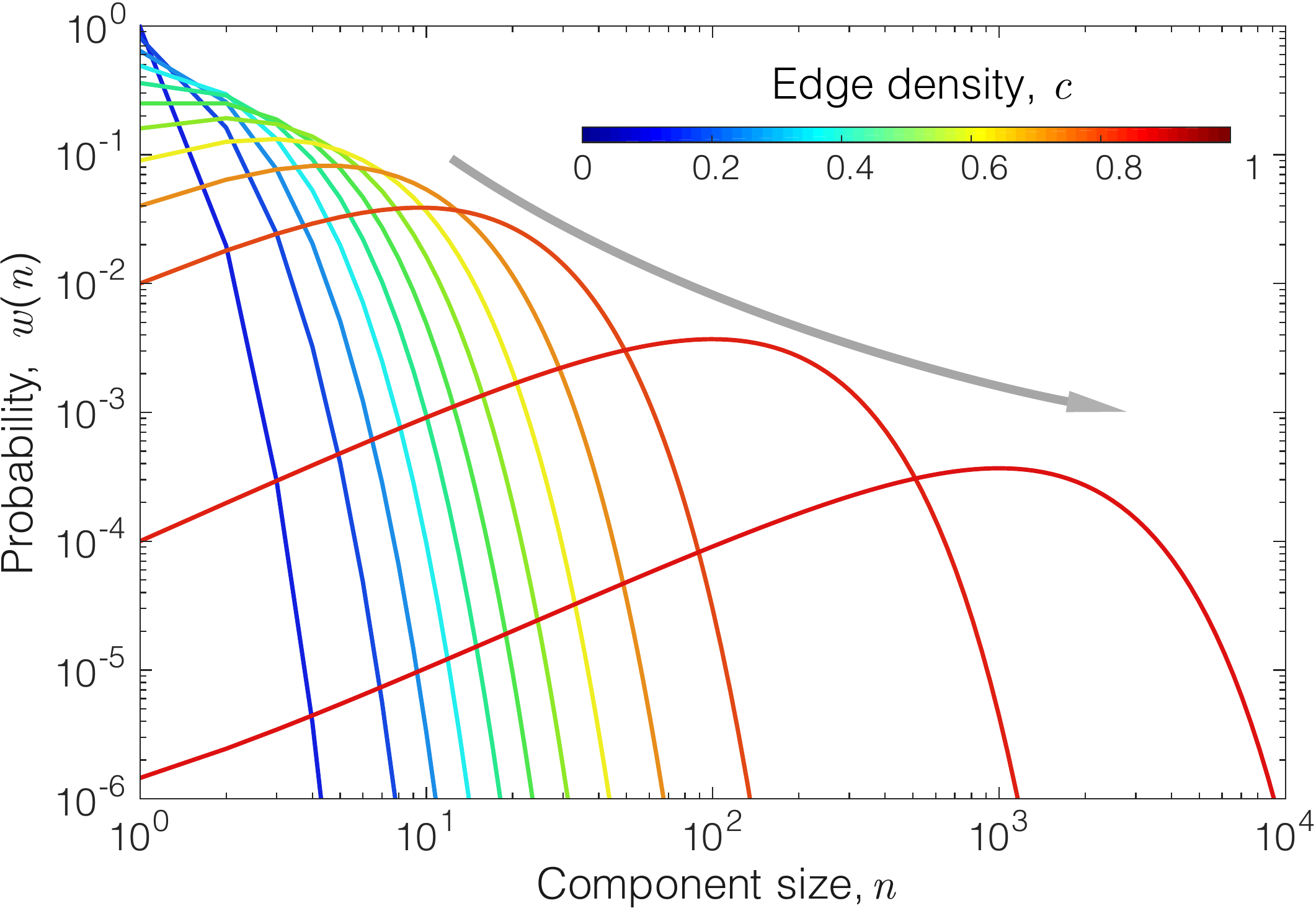} 
 \caption{Evolution of the size distribution of connected components for a system with $f_2=1$  and various values of the edge density.
The giant component emerges asymptotically at infinite time ($c \rightarrow 1$). }
\label{fig:single2F}
\end{center}
\end{figure}

\emph{Example 1} We consider vertices with at most degree 2, that is $f_2=1$. Graphs generated by such a process are always linear and, according to \eqref{eq:col1}, the giant component is reachable only asymptotically at $t\rightarrow \infty$. Furthermore, a small perturbation,  $f_1=\varepsilon,\; f_2=1-\varepsilon,$ prevents emergence of the giant component even at infinite time (see points A at barycentric plot of configurations, Figure~\ref{fig:bar}). The component-size distribution is illustrated in Figure~\ref{fig:single2F}. One may notice the constant ``drift'' (as indicated with an arrow) of the distribution towards larger values of components sizes. The distribution features the exponential asymptote at any $t>0$.

\begin{figure}[htbp]
\begin{center}
a. \includegraphics[width=0.7\textwidth]{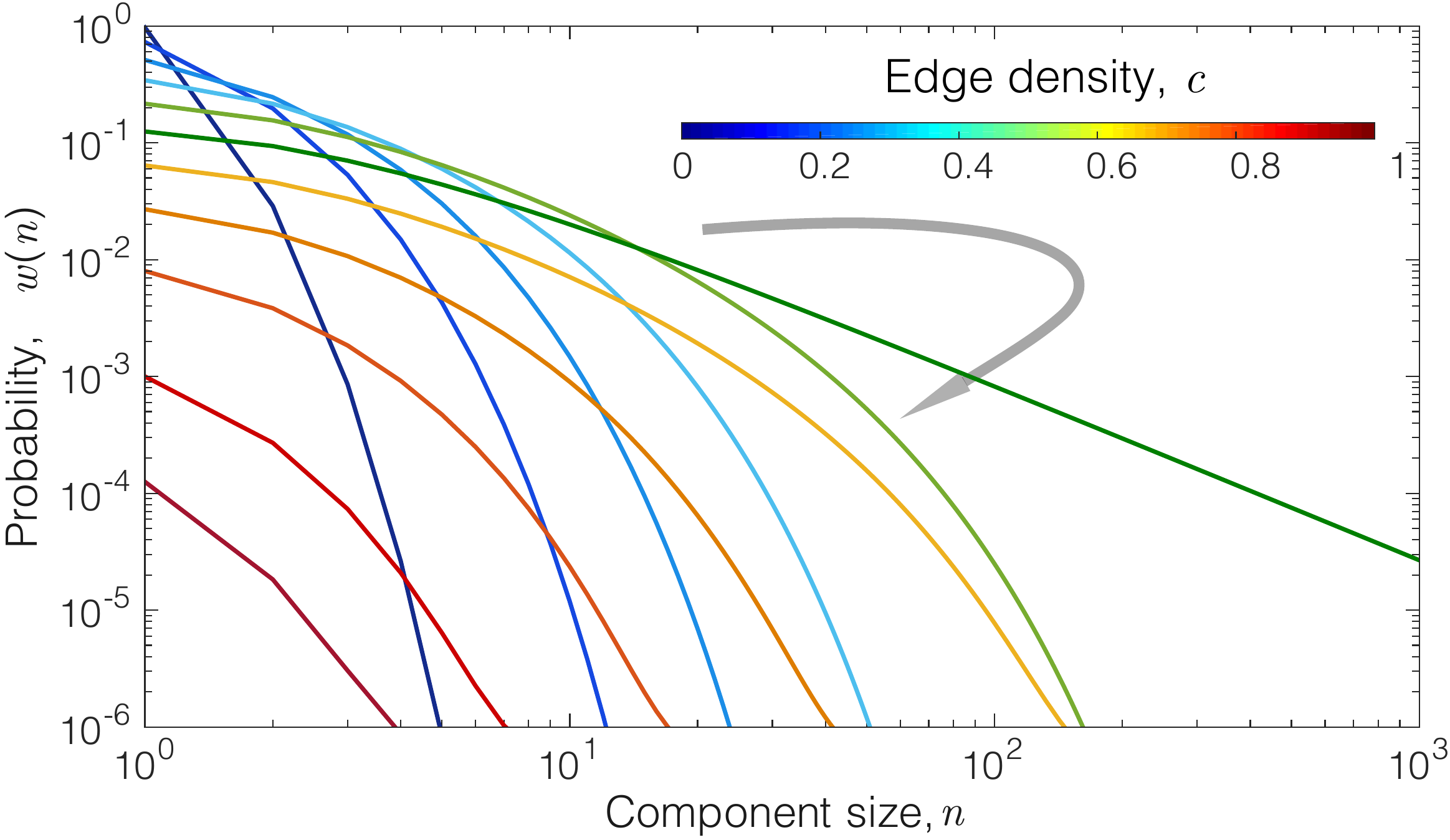} \\
b. \includegraphics[width=0.7\textwidth]{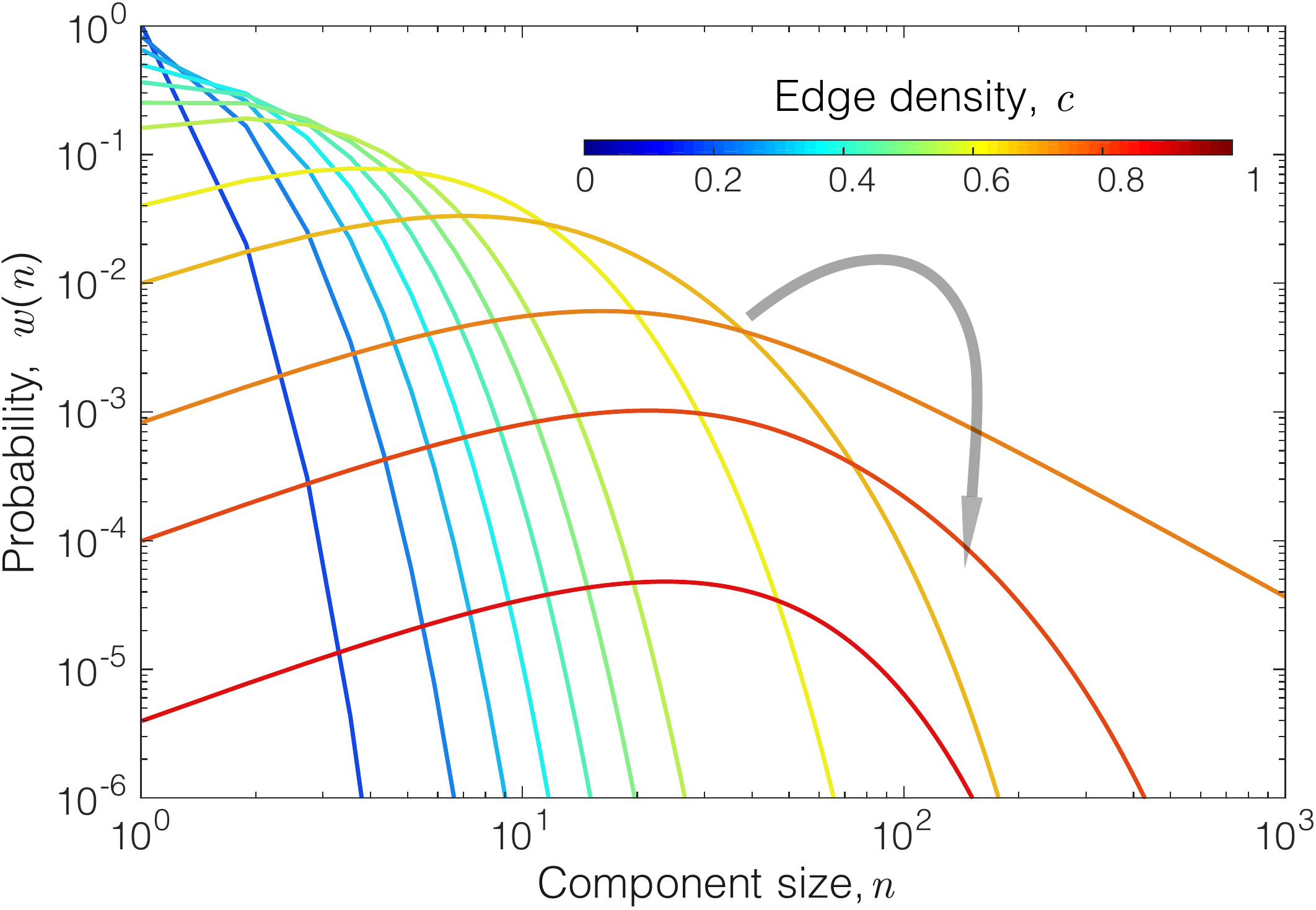} 
 \caption{Evolution of the size distribution of connected components for: \emph{a)}  a system with $f_3=1$, phase transition at $c_g=\frac{1}{2},$ and \emph{b)} a system with  $f_2=\frac{49}{50}, f_3=\frac{1}{50})$, phase transition at $c_g=\frac{101}{104}.$ 
  }
\label{fig:single3F}
\end{center}
\end{figure}
\begin{figure}[htbp]
\begin{center}
\includegraphics[width=\textwidth]{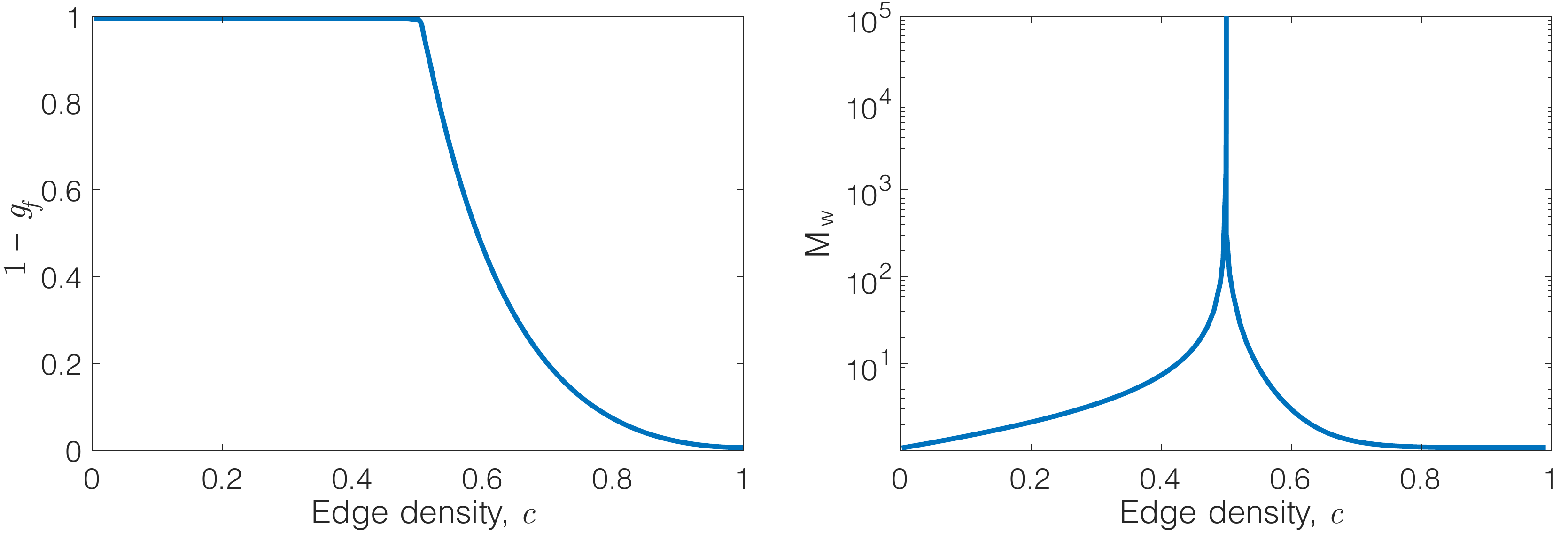} \\
 \caption{
 Emergence of the giant component in a system with $f_3=1$ that features phase transition at $c=\frac{1}{2}$.
 \emph{Left:} probability that a randomly sampled node belongs to a finite-size connected component. \emph{Right:} the expected size of connected components features a singularity at the phase transition.
  }
\label{fig:expected_size}
\end{center}
\end{figure}

\emph{Example 2} In this example we consider a system with $f_3=1.$ This random graph consists of three-functional vertices and features the phase transition at edge density $c=\frac{1}{2}$ (configuration C in Figure~\ref{fig:bar}).
The component-size distribution is illustrated in Figure~\ref{fig:single3F}a. Asymptotically, when $n\rightarrow \infty$, the component-size distribution switches between exponential decay ($0<c<0.5$), algebraic decay ($c=0.5$), and back to exponential decay again ($0.5<c<1$). Prior to the phase transition, the distribution 'drifts' to the right (expected component size becomes larger), and swings back to small expected component sizes at the end of the process. 
When edge density traverses the critical point $c_g=\frac{1}{2}$, the probability that a randomly sampled node belongs to finite-size component departs from one and the expected component size features a singularity, see Figure~\ref{fig:expected_size}.

As shown in Figure~\ref{fig:bar}, one may postpone the phase transition so it occurs anywhere between 0.5 and 1 by adding vertices of functionality 2 to the system. 
For instance, a mixture of vertices with functionalities two and three having fractions $f_2=\frac{49}{50}$ and $f_3=\frac{1}{50},$ as denoted by point B in Figure~\ref{fig:bar}, postpones the phase transition to $c_g=\frac{101}{104}\approx 0.97.$ The evolution of the size distribution for this case is depicted in Figure~\ref{fig:single3F}b.

While vertices of degree two postpone the phase transition, vertices of degree one may prevent it by ``consuming'' all available edges in a single connected component and thus locking its size finite. For this reason vertices of degree one are called \emph{termination agents} within the chemical context. Depending on what is the degree of the other species, the probability of randomly selecting a component may feature regular oscillations. For instance, in a dense, $c=1$, mixture of $m$-functional and one-functional vertices, connected components can take their sizes only from  
$$n\in \{2\}\cup \{k m-k+2\;|\;k=1,2,\dots \}.$$ 
Here we rely on the fact that non-giant components do not contain cycles\cite{newman2001}. Whereas when edge density $c<1,$ the sizes of connected components are not restricted to this set and, as is demonstrated in the next example, the transition of the size distribution from $c<1$ to $c=1$ is non-trivial. 
\begin{figure}[htbp]
\begin{center}
 a.
 \includegraphics[width=0.6\textwidth]{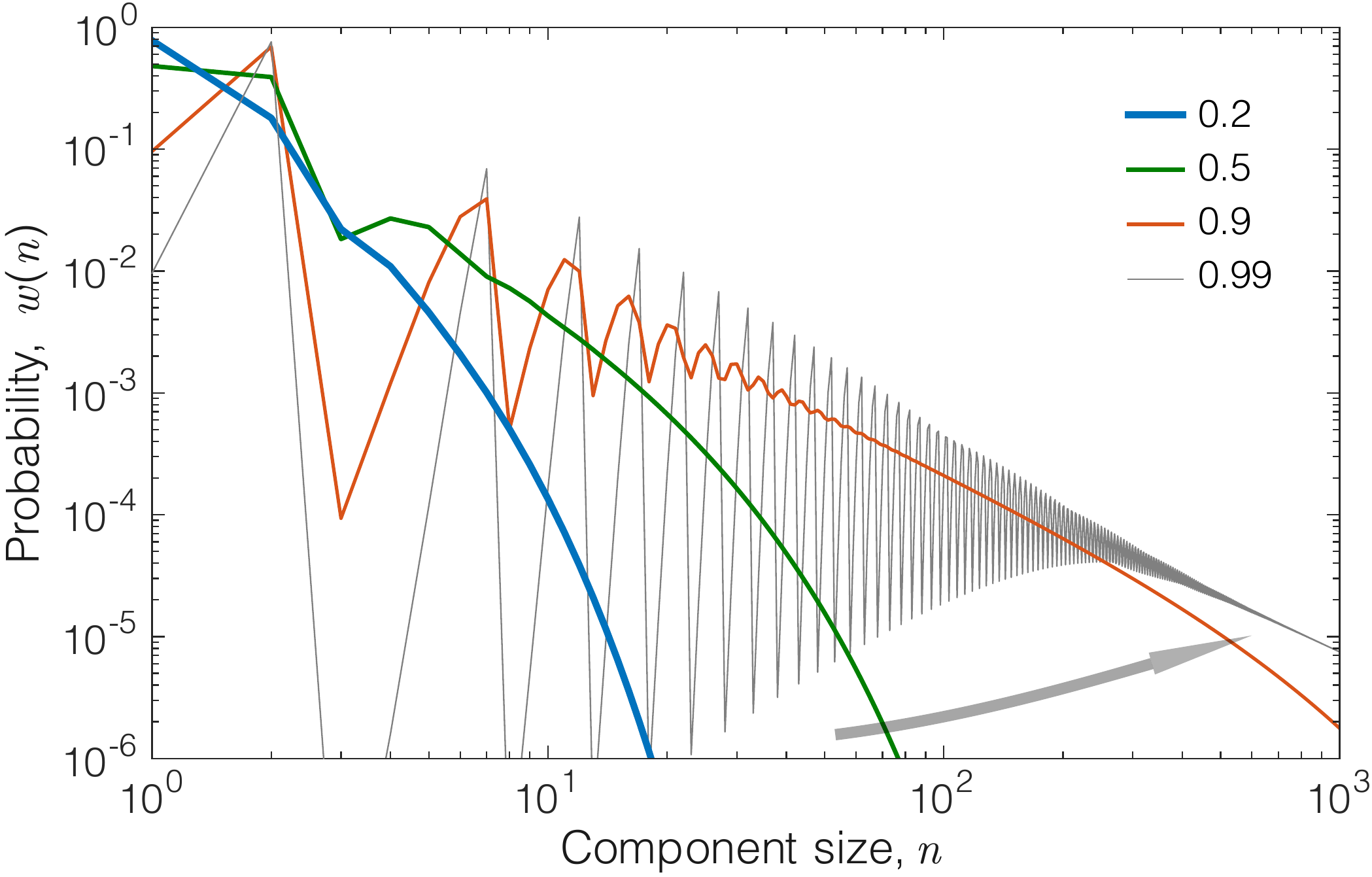} \\
 b. \includegraphics[width=0.6\textwidth]{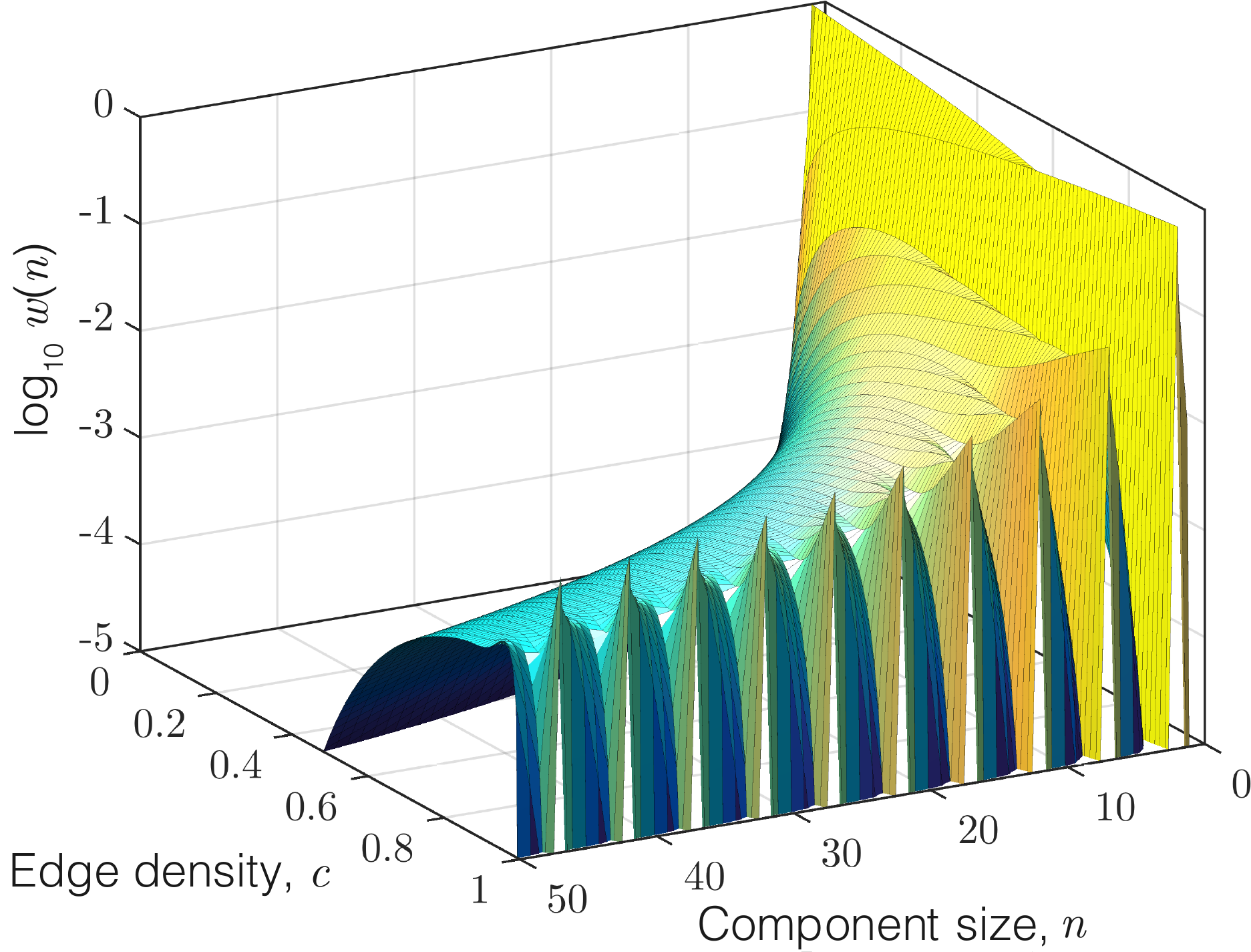} \\ 
 \caption{The size distributions of connected components for a system with $f_1=\frac{24}{25},\; f_6=\frac{1}{25}$ as predicted by the theory. \emph{a)} The size distributions at a few instances of time. \emph{b)} A surface representing the evolution of the size distribution during the whole time-continuous process, $c \in [0,1].$
   }
\label{fig:1and6}
\end{center}
\end{figure}

\begin{figure}[htbp]
\begin{center}
 \includegraphics[width=0.7\textwidth]{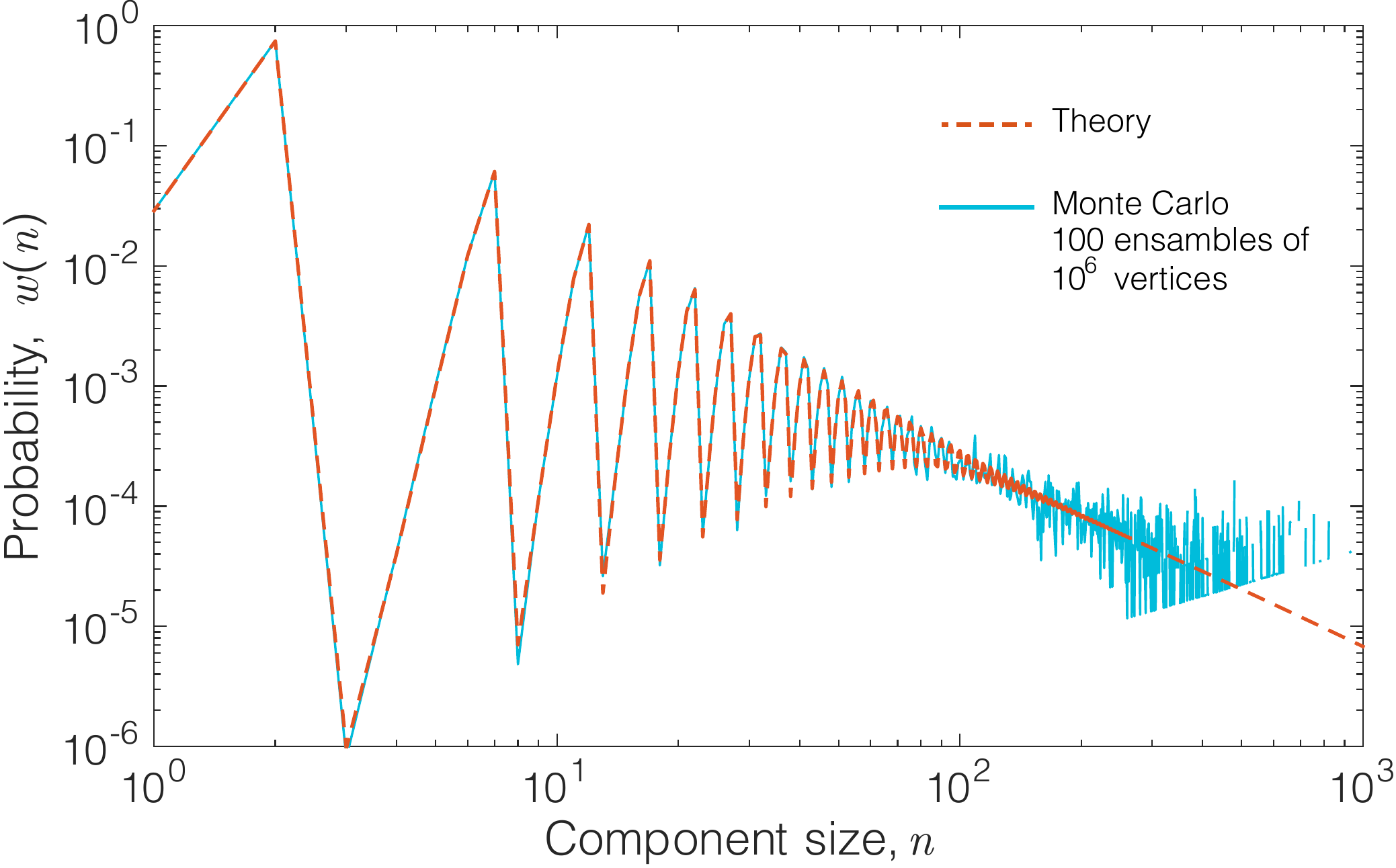} 
\caption{ The size distribution of connected components for a system with $f_1=\frac{24}{25},\; f_6=\frac{1}{25}$ at edge density $c=0.97$ is obtained with two different methods: (\emph{red line:})  the theory; (\emph{blue line:})  Monte Carlo simulations of a network with $10^6$ vertices. The simulation data is averaged over 100 simulation runs.
  }\label{fig:mc}
\end{center}
\end{figure}

\emph{Example 3} We consider a mixture of one- and six-functional  vertices present with fractions $f_1=\frac{24}{25},\;f_6=\frac{1}{25}$. This distribution of functionalities features the phase transition at $c=1.$ As illustrated in Figure~\ref{fig:1and6}a, the size distribution decays monotonically at low edge densities, but switches to oscillations as $c$ approaches 1. The switch itself is gradual as can be seen in Figure~\ref{fig:1and6}b. In Figure~\ref{fig:mc}, the theoretical results are compared to component-size distribution 
generated by Monte Carlo (MC) computations. The theory and MC data are in a perfect agreement; however, despite extensive size of MC computations (100 ensembles of size $10^6$ vertices), the MC resolution in the tail of the distributions remains poor.

\begin{acknowledgements}
This work is part of the project number 639.071.511, which is financed by the Netherlands Organisation for Scientific Research (NWO) VENI. Some of the results published in this paper were obtained during work under PAinT (Paint Alterations in Time) project as part of the NWO Science4Arts Program. 
\end{acknowledgements}

\bibliographystyle{spmpsci}  
\bibliography{literature}

\end{document}